\documentclass{article}
\usepackage{graphicx} 
\usepackage[utf8]{inputenc}
\usepackage[T1]{fontenc}
\usepackage{amsmath,amssymb,amsthm,mathtools}
\usepackage[margin=1in]{geometry}
\usepackage{hyperref}
\usepackage{relsize}

\theoremstyle{plain}
\newtheorem{theorem}{Theorem}[section]
\newtheorem{lemma}[theorem]{Lemma}
\newtheorem{proposition}[theorem]{Proposition}
\newtheorem*{conjecture}{Conjecture}

\theoremstyle{remark}
\newtheorem{remark}[theorem]{Remark}

\newcommand{\R}{\mathbb{R}}
\newcommand{\C}{\mathbb{C}}
\newcommand{\N}{\mathbb{N}}

\title{Note on the second derivatives of symbols for Gelfand--Shilov spaces}
\author{Héctor Ariza\footnote{E-mail: harirem@upvnet.upv.es}}

\date{November 2026}

\begin{document}

\maketitle

\begin{abstract}
    We investigate whether the continuity of the composition operator $C_\psi(f) = f \circ \psi$ on the Beurling-type Gelfand--Shilov space $\Sigma_d(\R)$ ($d > 1$) imposes a growth condition on the second derivative $\psi''$. In particular, we prove that $|\psi''(x)|$ must be bounded by some constant $C_\varepsilon$ times $(1+|x|)^{1-\frac{1}{d}+\varepsilon}$, for every $\varepsilon>0$. We provide an explicit symbol $\psi$ for $\Sigma_d(\mathbb{R})$ whose unbounded second derivative attains polynomial growth of order $1-\frac{1}{d}$, showing the quasi-optimality of the necessary condition. 
\end{abstract}\vspace{0.1cm}

\vspace{0.05cm}


\section{Introduction}
In recent years, there has been a growing interest in composition operators on the Schwartz class (see, for instance, \cite{akj,fgj,advances,gj}) and in the ultradifferentiable setting (see, for instance, \cite{albaneseariza, AFG24,WeightedModulationSpaces,Iterates,fg, rs}). Despite this, no full characterization of the continuity of the composition operator $C_\psi(f) = f \circ \psi$ acting on the Beurling-type Gelfand--Shilov space $\Sigma_d(\R)$ ($d > 1$) is known. In this note, we narrow this gap. In particular, we provide a sharp necessary condition on the second derivative of $\psi$. \\ 



For $d > 1$, the classical Beurling-type Gelfand--Shilov space $\Sigma_d(\R)$ is defined as the Fr\'echet space
\[
\Sigma_d(\R) := \left\{ f \in C^\infty(\R) : \|f\|_{h,\lambda,d} := \sup_{x \in \R,\, n \in \N_0} \frac{|f^{(n)}(x)|}{h^n (n!)^d} e^{\lambda|x|^{1/d}} < \infty \quad \forall h>0,\hspace{0.05cm}\forall\lambda > 0 \right\},
\]
equipped with the locally convex topology generated by the family of seminorms $\{\|\cdot\|_{h,\lambda,d}\}_{h>0, \hspace{0.02cm}\lambda>0}$. \\

Equivalently, the topology of $\Sigma_d(\R)$ is generated by the one-parameter family of seminorms (see \cite[Lemma 2.11]{asensio_jornet}, \cite{GelfandShilovBook}, \cite[Theorem 2.4]{paley} or \cite[Theorem 4.8]{seminormas}):
\begin{equation}\label{eq:seminorm}
\|f\|_{h,d}
:=
\sup_{y\in\mathbb{R},\;j,q\in\mathbb{N}_0}
\frac{(1+|y|)^q\,|f^{(j)}(y)|}{h^{j+q}\,(j+q)!^{\,d}}, \textit{ for each } h>0.
\end{equation}
For $d>1$, the space $\Sigma_d(\R)$ contains non-trivial compactly supported functions, is invariant under translation, differentiation, and the Fourier transform, and is a Montel space (that is, closed bounded subsets are compact). See \cite{albanesemele} and the references therein for a survey on Gelfand--Shilov spaces. \\

In \cite[Theorem 3.9]{AFG24}, it was shown that if the composition operator $C_\psi(f) = f \circ \psi$ maps $\Sigma_d(\R)$ continuously into itself, then $\psi'$ is bounded on $\R$. Since $\psi$ cannot be constant, from now on we set $M_1 := \sup_{y \in \R} |\psi'(y)|>0$ and $C_1 := |\psi(0)| + M_1 > 0$, so that
\begin{equation}\label{eq:psi_prime_bound}
|\psi(y)| \le |\psi(0)| + M_1|y| \le C_1(1 + |y|) \qquad \forall y \in \R.
\end{equation}

To the best of our knowledge, no explicit symbol for $\Sigma_d(\R)$ with unbounded second derivative has appeared in the literature. In Section~\ref{ExampleSection}, we construct a one-parameter family of such symbols whose simplicity allows sharp asymptotic estimates for all higher-order derivatives. \\

We are now interested in the behavior of $\psi''$. We will prove that, although $\psi''$ may be unbounded (Theorem \ref{thm:main}), it has polynomial growth of order at most $1-\frac{1}{d}+\varepsilon$ for every $\varepsilon>0$ (Theorem \ref{mainth}).

\section{Necessary condition on $\psi''$}\label{mainsec}
Our strategy builds on \cite[Theorem~3.9]{AFG24}. However, \cite[Lemma~2.1]{AFG24}, which allows one to collapse Fa\`a di Bruno's formula in that proof to a single term containing a power of $\psi'$, cannot be adapted to $\psi''$. For this reason, we need the following technical lemma:
\begin{lemma}\label{lem:linear_subtraction}
Suppose $C_\psi : \Sigma_d(\R) \to \Sigma_d(\R)$ is continuous, and let $M_1 = \sup_{y \in \R} |\psi'(y)| < \infty$. For any fixed $y_m \in \R$, define the following function:
\[
\tilde{\psi}_m(y) := \psi(y) - \psi'(y_m)(y - y_m) \hspace{0.15cm} \forall \hspace{0.03cm} y\in\mathbb{R}.
\]
Then $\tilde{\psi}_m(y_m) = \psi(y_m)$, $\tilde{\psi}_m'(y_m) = 0$, and $\tilde{\psi}_m'' = \psi''$. Moreover, for every bounded subset $\mathcal{B} \subset \Sigma_d(\R)$ and every $\lambda>0,\hspace{0.05cm}h > 0$, there exists a constant $C_{h,\lambda,M_1,\mathcal{B}} \geq 1$ independent of $m$ (of $y_m$) and $n$ such that
\begin{equation}\label{eq:flattened_bound}
\left|(g \circ \tilde{\psi}_m)^{(n)}(y_m)\right| \le C_{h,\lambda,M_1,\mathcal{B}}\hspace{0.1cm} (2h)^n \hspace{0.1cm}(n!)^d \hspace{0.1cm} e^{-\lambda|y_m|^{1/d}} \qquad \forall g \in \mathcal{B},\ \forall n \in \N_0.
\end{equation}
\end{lemma}

\begin{proof}
The derivative identities at $y = y_m$ follow immediately. Set $\xi_m := \psi'(y_m)$, so $|\xi_m| \le M_1$ and $\tilde{\psi}_m(y) = \psi(y) - \xi_m(y - y_m)$, for all $y\in\mathbb{R}$.\\

By the Fourier inversion formula, every function $g \in \Sigma_d(\R)$ satisfies
\[
g(\tilde{\psi}_m(y)) = \frac{1}{2\pi} \int_{\R} \hat{g}(\eta)\hspace{0.09cm} e^{-i\eta \xi_m (y - y_m)}\hspace{0.09cm} e^{i\eta \psi(y)} \, d\eta,
\] for all $y\in\mathbb{R}$.\\

Differentiating $n$ times with respect to $y$ at $y = y_m$ via Leibniz's rule (justified by the rapid decay of $\hat{g}$ and the smoothness of $\psi$) gives
\begin{align}
(g \circ \tilde{\psi}_m)^{(n)}(y_m) &= \sum_{j=0}^{n}\hspace{0.01cm} \binom{n}{j}\hspace{0.07cm} (-\xi_m)^{n-j} \hspace{0.07cm}\frac{1}{2\pi} \hspace{0.07cm}\int_{\R} \hspace{0.07cm}(i\eta)^{n-j} \hspace{0.07cm}\hat{g}(\eta) \left.\hspace{0.07cm} \frac{d^j}{dy^j}\left( e^{i\eta \psi(y)} \right) \right|_{y=y_m} d\eta \nonumber \\
&= \sum_{j=0}^{n}\hspace{0.07cm} \binom{n}{j}\hspace{0.09cm} (-\xi_m)^{n-j}\hspace{0.09cm} \bigl( C_\psi(D^{n-j} g) \bigr)^{(j)}(y_m), \label{eq:leibniz_fourier}
\end{align} where $D=\frac{d}{dx}$. The second equality above holds because of the following identity:
\[
\bigl( C_\psi(D^{n-j} g) \bigr)(y) = (D^{n-j} g)(\psi(y)) = \frac{1}{2\pi} \hspace{0.08cm}\int_{\R} \hspace{0.06cm}(i\eta)^{n-j} \hspace{0.1cm}\hat{g}(\eta) \hspace{0.1cm}e^{i \eta \psi(y)} \, \hspace{0.08cm}d\eta,
\] since $\widehat{D^{n-j} g}(\eta) = (i\eta)^{n-j}\hat{g}(\eta)$. Therefore, taking the $j$-th derivative with respect to $y$ on both sides (interchanging differentiation and integration justified by the rapid decay of $\hat{g}$) and evaluating at $y = y_m$ yields
\[
\bigl( C_\psi(D^{n-j} g) \bigr)^{(j)}(y_m) = \frac{1}{2\pi} \int_{\R} (i\eta)^{n-j} \hat{g}(\eta) \left. \frac{d^j}{dy^j}\left( e^{i\eta \psi(y)} \right) \right|_{y=y_m} d\eta.
\]

Let $\lambda>0$. Fix $h > 0$ and define $R\equiv R(h) := \frac{h}{M_1 + 1} > 0$. Consider the family of continuous linear operators $\{T_p\}_{p \in \N_0}$ on $\Sigma_d(\R)$ given by
\[
T_p(g) := \frac{1}{R^p \hspace{0.07cm}(p!)^d}\hspace{0.07cm} C_\psi(D^p g).
\] We claim that $\{T_p(g)\}_{p \ge 0}$ is bounded in $\Sigma_d(\R)$ for every $g \in \Sigma_d(\R)$. Indeed, for any $\rho_1, \rho_2 > 0$, set $\tilde{\rho} := \min(\rho_1, R)/2^d$. Since $\tilde{\rho} \le R/2^d$, $\tilde{\rho} \le \rho_1/2^d$, and $(p+k)! \le 2^{p+k} \hspace{0.07cm}p! \hspace{0.07cm} k!$, we obtain
\[
\frac{|(D^p g)^{(k)}(x)|}{R^p \hspace{0.07cm}(p!)^d \hspace{0.07cm}\rho_1^k \hspace{0.07cm}(k!)^d} \hspace{0.09cm}e^{\rho_2 |x|^{1/d}} \le \frac{|g^{(p+k)}(x)|}{\tilde{\rho}^{p+k} \hspace{0.07cm}((p+k)!)^d} \hspace{0.09cm}e^{\rho_2 |x|^{1/d}} \le \|g\|_{\tilde{\rho}, \rho_2, d} < \infty,
\] for all $p, k \in \N_0$ and $x \in \R$. Thus $\left\{ \frac{1}{R^p\hspace{0.07cm} (p!)^d}\hspace{0.07cm} D^p g \right\}_{p \ge 0}$ is bounded in $\Sigma_d(\R)$, and by continuity of $C_\psi$, $\{T_p(g)\}_{p \ge 0}$ is bounded in $\Sigma_d(\R)$, for every $g\in \Sigma_d(\R)$. Since $\Sigma_d(\R)$ is a Fr\'echet space, the Banach--Steinhaus theorem implies that $\{T_p\}_{p \ge 0}$ is equicontinuous. Consequently, for any bounded set $\mathcal{B} \subset \Sigma_d(\R)$, the set $\bigcup_{p \ge 0} T_p(\mathcal{B})$ is bounded in $\Sigma_d(\R)$.\\

In particular, there exists $K = K_{h,\lambda,M_1,\mathcal{B}} > 0$ such that $\|T_p(g)\|_{h,\lambda,d} \le K$ for all $p \ge 0$ and $g \in \mathcal{B}$. Using this bound with $p = n-j$ for the $j$-th derivative at $y_m$ gives
\[
\left| \bigl( C_\psi(D^{n-j} g) \bigr)^{(j)}(y_m) \right| \le K \hspace{0.09cm}R^{n-j} \hspace{0.09cm}((n-j)!)^d\hspace{0.09cm} h^j \hspace{0.09cm}(j!)^d\hspace{0.09cm} e^{-\lambda|y_m|^{1/d}}.
\]
Inserting this bound into \eqref{eq:leibniz_fourier} and using $((n-j)!)^d \hspace{0.07cm}(j!)^d \le (n!)^d$ and $M_1 R + h = \frac{M_1 h}{M_1 + 1} + h < 2h$, we obtain
\begin{align*}
|(g \circ \tilde{\psi}_m)^{(n)}(y_m)| &\le K\hspace{0.09cm} (n!)^d \hspace{0.09cm}e^{-\lambda|y_m|^{1/d}} \hspace{0.09cm}\sum_{j=0}^{n} \binom{n}{j} \hspace{0.09cm}(M_1 R)^{n-j} \hspace{0.09cm}h^j \\
&= K \hspace{0.09cm}(M_1 R + h)^n \hspace{0.09cm}(n!)^d\hspace{0.09cm} e^{-\lambda|y_m|^{1/d}} \\
&\le K \hspace{0.09cm}(2h)^n \hspace{0.09cm}(n!)^d \hspace{0.09cm}e^{-\lambda |y_m|^{1/d}},
\end{align*}
establishing \eqref{eq:flattened_bound} with $C_{h,\lambda,M_1,\mathcal{B}}=\max\{1,K\}$.
\end{proof}

We are now ready to prove the main result of this note:
\begin{theorem} \label{mainth}
Let $d > 1$ and suppose $C_\psi : \Sigma_d(\R) \to \Sigma_d(\R)$ is well-defined and continuous. Then for every $\varepsilon > 0$, there exists a constant $C_\varepsilon > 0$ such that
\begin{equation}\label{objetivoEpsilon}
|\psi''(y)| \le C_\varepsilon (1 + |y|)^{1 - \frac{1}{d} + \varepsilon} \qquad \forall y \in \R.
\end{equation}
\end{theorem}

\begin{proof}
Fix $\varepsilon>0$. We assume without loss of generality that $0<\varepsilon<1$. Consider the Fr\'echet sequence space
\[
\Lambda_d := \left\{ a = (a_\ell)_{\ell \in \N_0} \subset \C : \|a\|_H := \sup_{\ell \in \N_0} \frac{|a_\ell|}{H^\ell (\ell!)^d} < \infty \quad \forall H > 0 \right\}.
\]
Let $r\in\mathbb{N}$, and define the sequence $\tilde{b}^{(r)} = (\tilde{b}_\ell^{(r)})_{\ell \in \N_0}$ by
\[
\tilde{b}_\ell^{(r)} := \begin{cases}
(r!)^{d - \varepsilon/2} & \text{if } \ell = r, \\
0 & \text{if } \ell \neq r.
\end{cases}
\]
For every $H > 0$, we have
\[
\sup_{r \ge 1} \|\tilde{b}^{(r)}\|_H = \sup_{r \ge 1} \frac{(r!)^{d - \varepsilon/2}}{H^r (r!)^d} = \sup_{r \ge 1} \frac{1}{H^r (r!)^{\varepsilon/2}} =: K_{H,\varepsilon} < \infty,
\]
so $\{\tilde{b}^{(r)}\}_{r \ge 1}$ is a bounded subset of $\Lambda_d$. Since $d > 1$, Petzsche's theorem (see \cite{petzsche} or \cite{AFG24} and the references therein) implies that the Borel map $B : \Sigma_d(\R) \to \Lambda_d$, defined by $$B(\tilde{\phi}) := \left(\tilde{\phi}^{(\ell)}(0)\right)_{\ell \in \N_0},$$ is a continuous linear surjective function between Fr\'echet spaces. By \cite[26.22]{meisevogt} and the fact that $\Lambda_d$ is Montel, there is a bounded sequence $\{\tilde{\phi}_r\}_{r=1}^\infty \subset \Sigma_d(\R)$ such that $B( \tilde{\phi}_r)=\tilde{b}^{(r)}$, for all $r\in\mathbb{N}$, and hence, $\tilde{\phi}_r^{(\ell)}(0) = \tilde{b}_\ell^{(r)}$ for all $\ell \in \N_0$, $r\in\mathbb{N}$. In particular, for every $\lambda > 0$ and $h>0$,
\begin{equation}\label{eq:lifted_bound}
M_{h,\lambda, d} := \sup_{r \ge 1} \|\tilde{\phi}_r\|_{h,\lambda, d} < \infty.
\end{equation}

Let $y_m \in \R$ be an arbitrary number with $|y_m| \ge 3^d$, set $x_m := \psi(y_m)$, and define the integer index
\[
r = r(m) := \lfloor |y_m|^{1/d} \rfloor \ge 3.
\]
Set $\lambda_m := \ln \ln r > 0$, and define the translated and rescaled test function
\[
g_m(x) := e^{-\lambda_m |x_m|^{1/d}}\hspace{0.09cm} \tilde{\phi}_r(x - x_m) \hspace{0.17cm} \forall\hspace{0.05cm} x\in\mathbb{R}.
\]
We claim that $\mathcal{B} := \{g_m : |y_m| \ge 3^d\}$ is a bounded subset of $\Sigma_d(\R)$. It suffices to verify that every sequence $\{g_{m_k}\}_{k \ge 1} \subset \mathcal{B}$ has a bounded subsequence in $\Sigma_d(\R)$. Notice that the sequence $\{\tilde{\phi}_{r(m_k)}\}_{k \ge 1}$ is bounded in $\Sigma_d(\R)$ by \eqref{eq:lifted_bound}. Passing to a subsequence if necessary, we distinguish two cases:
\begin{itemize}
    \item If $|y_{m_k}| \to \infty$ as $k \to \infty$, then $r(m_k) \to \infty$ and $\lambda_{m_k} = \ln \ln r(m_k) \to +\infty$.  
     Since $d > 1$, the function $t \mapsto t^{1/d}$ is subadditive on $[0, \infty[$. Therefore, $$\vert{}x\vert{}^{1/d} \le \vert{}x - x_{m_k}\vert{}^{1/d} + \vert{}x_{m_k}\vert{}^{1/d} \qquad \forall x \in \mathbb{R},\ \forall k \ge 1.$$ Fix $\lambda > 0$. Then $$\lambda \vert{}x\vert{}^{1/d} - \lambda_{m_k}\vert{}x_{m_k}\vert{}^{1/d} \le \lambda \vert{}x - x_{m_k}\vert{}^{1/d} + (\lambda - \lambda_{m_k})\vert{}x_{m_k}\vert{}^{1/d} \qquad \forall x \in \mathbb{R},\ \forall k \ge 1.$$ Obviously, $g_{m_k}^{(n)}(x) = e^{-\lambda_{m_k}\vert{}x_{m_k}\vert{}^{1/d}} \tilde{\phi}_{r(m_k)}^{(n)}(x - x_{m_k})$ for all $n \in \mathbb{N}_0$ and $x \in \mathbb{R}$. Fix $h>0$. We obtain for every $k \ge 1$:
\begin{equation*}
\begin{split}
\|g_{m_k}\|_{h,\lambda, d}
&= \sup_{\substack{x \in \mathbb{R}, \ n \in \mathbb{N}_0}} \frac{|g_{m_k}^{(n)}(x)|}{h^n \hspace{0.07cm}(n!)^d}\hspace{0.07cm} \hspace{0.07cm}e^{\lambda |x|^{1/d}} \
\\&= \sup_{\substack{x \in \mathbb{R}, \ n \in \mathbb{N}_0}} \frac{|\tilde{\phi}_{r(m_k)}^{(n)}(x - x_{m_k})|}{h^n \hspace{0.07cm}(n!)^d}\hspace{0.07cm} e^{\lambda |x|^{1/d} - \lambda_{m_k}|x_{m_k}|^{1/d}} \
\\&\le e^{(\lambda - \lambda_{m_k})\hspace{0.02cm}|x_{m_k}|^{1/d}} \hspace{0.07cm}\sup_{\substack{x \in \mathbb{R}, \ n \in \mathbb{N}_0}} \frac{|\tilde{\phi}_{r(m_k)}^{(n)}(x - x_{m_k})|}{h^n \hspace{0.07cm}(n!)^d} \hspace{0.07cm}e^{\lambda |x - x_{m_k}|^{1/d}} \
\\&= e^{(\lambda - \lambda_{m_k})|x_{m_k}|^{1/d}}\hspace{0.07cm} \|\tilde{\phi}_{r(m_k)}\|_{h,\lambda, d} \
\\ &\le e^{(\lambda - \lambda_{m_k})\hspace{0.02cm}|x_{m_k}|^{1/d}}\hspace{0.07cm} M_{h,\lambda, d}.
\end{split}
\end{equation*}

Because $\lambda_{m_k} \to +\infty$ as $k \to \infty$, there exists $k_0 \in \mathbb{N}$ such that $\lambda_{m_k} \ge \lambda$ for all $k \ge k_0$. For every $k \ge k_0$, we have $\lambda - \lambda_{m_k} \le 0$, so $e^{(\lambda - \lambda_{m_k})\vert{}x_{m_k}\vert{}^{1/d}} \le 1$ and hence $\Vert{}g_{m_k}\Vert{}_{h,\lambda, d} \le M_{h,\lambda, d}$. Defining$$K_\lambda := \max\Bigl\{1, \max \bigl\{e^{(\lambda - \lambda_{m_k})\vert{}x_{m_k}\vert{}^{1/d}} : 1 \le k \leq k_0\bigr\}\Bigr\} < \infty,$$ we conclude that$$\sup_{k \ge 1} \Vert{}g_{m_k}\Vert{}_{h,\lambda, d} \le K_\lambda \hspace{0.07cm}M_{h,\lambda, d} < \infty.$$ Since $\lambda > 0$ and $h>0$ were arbitrary, $\{g_{m_k}\}_{k \ge 1}$ is bounded in $\Sigma_d(\mathbb{R})$.\footnote{By \cite[Lemma~2.3]{AFG24}, $\{g_{m_k}\}_{k\ge 1}$ is bounded in $\Sigma_d(\R)$; we include the proof here for completeness.}
    \item If $\{|y_{m_k}|\}_{k\ge 1}$ is bounded, say $|y_{m_k}| \le Y$, then \eqref{eq:psi_prime_bound} implies that $$|x_{m_k}| = |\psi(y_{m_k})| \le C_1(1+Y) =: X < \infty,$$ for all $k\in\mathbb{N}$. Since $\lambda_{m_k} > 0$, subadditivity of $t \mapsto t^{1/d}$ on $[0,\infty[$ implies that
\[
\lambda |x|^{1/d} - \lambda_{m_k} |x_{m_k}|^{1/d} \le \lambda |x - x_{m_k}|^{1/d} + \lambda X^{1/d}
\]
for every $\lambda > 0$ and $x \in \mathbb{R}$. Clearly, $g_{m_k}(x) = e^{-\lambda_{m_k}\hspace{0.09cm}|x_{m_k}|^{1/d}} \hspace{0.09cm}\tilde{\phi}_{r(m_k)}(x - x_{m_k})$, for all $k\in\mathbb{N}$ and $x\in\mathbb{R}$. For every $h>0,\hspace{0.03cm}\lambda > 0$ we obtain
\begin{align*}
\|g_{m_k}\|_{h,\lambda, d} 
&= \sup_{\substack{x \in \mathbb{R} \\ n \in \mathbb{N}_0}} \frac{|g_{m_k}^{(n)}(x)|}{h^n\hspace{0.09cm} (n!)^d} \hspace{0.09cm}e^{\lambda |x|^{1/d}} \\
&= \sup_{\substack{x \in \mathbb{R} \\ n \in \mathbb{N}_0}} \frac{|\tilde{\phi}_{r(m_k)}^{(n)}(x - x_{m_k})|}{h^n\hspace{0.09cm} (n!)^d} \hspace{0.09cm}e^{\lambda |x|^{1/d} - \lambda_{m_k}|x_{m_k}|^{1/d}} \\
&\le e^{\lambda X^{1/d}} \hspace{0.09cm}\sup_{\substack{x \in \mathbb{R} \\ n \in \mathbb{N}_0}} \frac{|\tilde{\phi}_{r(m_k)}^{(n)}(x - x_{m_k})|}{h^n \hspace{0.09cm}(n!)^d}\hspace{0.09cm} e^{\lambda |x - x_{m_k}|^{1/d}} \\
&= e^{\lambda X^{1/d}} \|\tilde{\phi}_{r(m_k)}\|_{h,\lambda, d} \le e^{\lambda X^{1/d}} M_{h,\lambda, d},
\end{align*}
since $\sup_{k \ge 1} \|\tilde{\phi}_{r(m_k)}\|_{h,\lambda, d}\leq M_{h,\lambda, d}   < \infty$ (see \eqref{eq:lifted_bound}).

\end{itemize}

Consequently, in both cases the sequence $\{g_{m_k}\}_{k \ge 1} \subset \mathcal{B}$ has a bounded subsequence in $\Sigma_d(\R)$. Therefore, the set $\mathcal{B}$ is bounded in $\Sigma_d(\R)$. \\

Define $\tilde{\psi}_m$ as in Lemma \ref{lem:linear_subtraction} for each $y_m\in\mathbb{R}$. Now, we apply Fa\`a di Bruno's formula to compute $(g_m \circ \tilde{\psi}_m)^{(2r)}$ evaluated at $y_m$:
\begin{equation}\label{FaaDiBrunoPsiGorro}
    (g_m \circ \tilde{\psi}_m)^{(2r)}(y_m) = \sum_{k=1}^{2r} g_m^{(k)}(x_m) \sum_{\substack{k_1 + \dots + k_{2r} = k \\ k_1 + 2k_2 + \dots + 2r k_{2r} = 2r}} \frac{(2r)!}{k_1! \cdots k_{2r}!} \prod_{j=1}^{2r} \left( \frac{\tilde{\psi}_m^{(j)}(y_m)}{j!} \right)^{k_j}.
\end{equation}
Since $g_m^{(k)}(x_m) = e^{-\lambda_m |x_m|^{1/d}} \hspace{0.09cm}\tilde{\phi}_r^{(k)}(0) = 0$ for all $k \neq r$, the outer sum collapses to $k = r$. Moreover, since $\tilde{\psi}_m'(y_m) = 0$, every term in the inner sum corresponding to $(k_1,k_2,...,k_{2\hspace{0.01cm} r}) \in \mathbb{N}_0^{2r}$ satisfying $k_1+k_2+...+k_{2\hspace{0.01cm}r}=r$, $k_1+2k_2+...+2\hspace{0.02cm}r \hspace{0.02cm}k_{2\hspace{0.01cm}r}=2\hspace{0.02cm}r$ and $k_1\geq 1$ vanishes. The remaining constraints $k_2 + k_3 + \dots + k_{2r} = r$ and $2k_2 + 3k_3 + \dots + 2r k_{2r} = 2r$ imply $k_3 + 2k_4 + \dots + (2r-2)k_{2r} = 0$, and hence $k_2 = r$ and $k_3 = \dots = k_{2r} = 0$, since $k_3, ..., k_{2\hspace{0.02cm}r}\geq 0$. Therefore, only a single term remains in the sum in \eqref{FaaDiBrunoPsiGorro}:
\[
|(g_m \circ \tilde{\psi}_m)^{(2r)}(y_m)| = |g_m^{(r)}(x_m)| \hspace{0.09cm}\frac{(2r)!}{r!\, 2^r} \hspace{0.09cm}|\tilde{\psi}_m''(y_m)|^r = e^{-\lambda_m |x_m|^{1/d}}\hspace{0.09cm} (r!)^{d - \varepsilon/2} \hspace{0.09cm}\frac{(2r)!}{r!\, 2^r} \hspace{0.09cm}|\psi''(y_m)|^r,
\] for all $r=r(y_m)\in\mathbb{N}$ and $|y_m|\geq 3^d$. On the other hand, applying the upper bound \eqref{eq:flattened_bound} from Lemma~\ref{lem:linear_subtraction} at $n = 2r$ and $h=\lambda = 1$ yields
\[
|(g_m \circ \tilde{\psi}_m)^{(2r)}(y_m)| \le C_{1,1,M_1,\mathcal{B}}\hspace{0.1cm} 4^{r}\hspace{0.1cm} ((2r)!)^d\hspace{0.1cm} e^{-|y_m|^{1/d}}.
\]
Combining this upper bound with the preceding identity yields
\begin{equation}\label{psisegderupper}
|\psi''(y_m)|^r \le C_{1,1,M_1,\mathcal{B}}\hspace{0.09cm} 8^r \hspace{0.09cm}\frac{((2r)!)^{d-1}}{(r!)^{d - 1 - \varepsilon/2}} \hspace{0.09cm}\exp\left( \lambda_m |x_m|^{1/d} - |y_m|^{1/d} \right).
\end{equation}
For $|y_m| \ge 3^d$, \eqref{eq:psi_prime_bound} implies that $|x_m| = |\psi(y_m)| \le |\psi(0)| + M_1|y_m| \le C_1|y_m|$, where $C_1 = |\psi(0)| + M_1$. Furthermore, since $r = \lfloor |y_m|^{1/d} \rfloor \ge 3$, we have $|y_m|^{1/d} < r + 1 \le \frac{4}{3} r$, and hence,
\[
\frac{|x_m|^{1/d}}{r} \le \frac{C_1^{1/d}\hspace{0.08cm} |y_m|^{1/d}}{r} \le \frac{4}{3}\hspace{0.08cm} C_1^{1/d} =: K_0.
\]
Therefore, the following inequality holds:
\begin{equation*}
    \begin{split}
        \left[ \exp\left( \lambda_m |x_m|^{1/d} - |y_m|^{1/d} \right) \right]^{1/r} &= \exp\left( \frac{(\ln \ln r)\hspace{0.06cm} |x_m|^{1/d}}{r} - \frac{|y_m|^{1/d}}{r} \right) \\&\le e^{K_0 \ln \ln r - 1} \le (\ln r)^{K_0}.
    \end{split}
\end{equation*}
We observe that $\frac{(2r)!}{(r!)^2}=\binom{2r}{r}\le 4^r$ and also $(r!)^{1/r} \le r$, 
for all $r\in\mathbb{N}$. \\

Taking the $r$-th root of both sides of the inequality \eqref{psisegderupper}, using the estimates above and the fact that $(\ln r)^{K_0} \le A_{K_0, \varepsilon}\hspace{0.09cm} r^{\varepsilon/2}$ for all $r \ge 3$, where $$A_{K_0, \varepsilon} := \sup_{t \ge 3} \left((\ln t)^{K_0} t^{-\varepsilon/2}\right) < \infty,$$ we deduce that
\begin{align*}
|\psi''(y_m)| &\le 8\hspace{0.09cm}C_{1,1,M_1,\mathcal{B}}^{1/r}\hspace{0.09cm}  \left( \left( \frac{(2r)!}{(r!)^2} \right)^{d-1} \cdot (r!)^{d - 1 + \varepsilon/2} \right)^{1/r} \hspace{0.06cm}(\ln r)^{K_0} \\
&\le 8\hspace{0.03cm}\cdot\hspace{0.03cm} 4^{d-1} \hspace{0.1cm}C_{1,1,M_1,\mathcal{B}}^{1/r}\hspace{0.1cm} \bigl((r!)^{1/r}\bigr)^{d - 1 + \varepsilon/2} \hspace{0.09cm}A_{K_0, \varepsilon} \hspace{0.1cm}r^{\varepsilon/2} \\
&\le \tilde{C}_\varepsilon \hspace{0.1cm}r^{d - 1 + \varepsilon},
\end{align*} where $\tilde{C}_\varepsilon=8\hspace{0.03cm}\cdot\hspace{0.03cm} 4^{d-1} \hspace{0.1cm}C_{1,1,M_1,\mathcal{B}}\hspace{0.09cm}A_{K_0, \varepsilon}.$
Since $r = \lfloor |y_m|^{1/d} \rfloor \le |y_m|^{1/d}$, we obtain 
\[
|\psi''(y_m)| \le \tilde{C}_\varepsilon \hspace{0.07cm}|y_m|^{1 - \frac{1}{d} + \frac{\varepsilon}{d}} \le \tilde{C}_\varepsilon\hspace{0.07cm} (1 + |y_m|)^{1 - \frac{1}{d} + \varepsilon},
\]
as $d > 1$ implies $\frac{\varepsilon}{d} < \varepsilon$. Because $|y_m| \ge 3^d$ was arbitrary and $\psi''$ is continuous on the compact interval $[-3^d, 3^d]$, \eqref{objetivoEpsilon} holds for all $y \in \R$ with
\[
C_\varepsilon := \max\left( \tilde{C}_\varepsilon, \sup_{|y| \le 3^d} |\psi''(y)| \right).
\]
\end{proof}

\section{Quasi-optimality of the polynomial growth upper bound for $\psi''$}\label{ExampleSection}
In this section we provide an explicit example of a symbol $\psi$ for $\Sigma_d(\mathbb{R})$ whose second derivative attains the upper bound set in Theorem \ref{mainth} (with $\varepsilon=0$).\\ 

Let $d>1$ and $1 < \alpha \le 2 - \frac{1}{d}$. Set
$\beta := \alpha - 1 \in \left]0,\, 1-\frac{1}{d}\right]$ and define $\phi(x) := (1+x^2)^{\alpha/2}$, $g(x) := \sin(\phi(x))$, for $x\in\mathbb{R}$, and $\psi \in C^\infty(\mathbb{R})$ by
\begin{equation}\label{counterexamplePsi}
\psi(x) := \int_0^x \bigl(1 + g(t)\bigr)\,dt, \qquad \forall \hspace{0.05cm} x\in\mathbb{R}.
\end{equation}

Clearly, $\psi'(x) = 1 + \sin\bigl((1+x^2)^{\alpha/2}\bigr)$, for all $x\in\mathbb{R}$, and $\psi(0)=0.$ \\

For our purposes in this section, it is convenient to work with the equivalent family of seminorms for $\Sigma_d(\mathbb{R})$ given in \eqref{eq:seminorm}.


\begin{theorem}\label{thm:main}
Let $d>1$ and $1<\alpha\le 2-\frac{1}{d}$. If $\psi$ is the function given in \eqref{counterexamplePsi}, then 
\begin{enumerate}
\item $\psi'(x) \in [0,2]$ for all $x\in\mathbb{R}$, whereas $\psi''$ is unbounded.\\ In fact, $\psi''(x_j)=\alpha\hspace{0.07cm}x_j\hspace{0.01cm}\left(1+x_j^2\right)^{(\alpha-2)/2}\sim \alpha \hspace{0.07cm}x_j^{\alpha-1}$ as $j\to\infty$, for $x_j = \left((2\pi j)^{2/\alpha} - 1\right)^{1/2},$ with $j\in\mathbb{N}$.
\item The composition operator $C_\psi: \Sigma_d(\mathbb{R})\to \Sigma_d(\mathbb{R})$ is continuous.
\end{enumerate}
\end{theorem}

\vspace{0.2cm}

We first establish some preliminary estimates:
\begin{proposition}\label{prop:growth}
The function $\psi$ given in \eqref{counterexamplePsi} belongs to $C^\infty(\mathbb{R})$, satisfies $0 \le \psi'(x) \le 2$ for all $x\in\mathbb{R}$ and
\begin{equation}\label{eq:comparable}
1 + |x| \le C_0\bigl(1 + |\psi(x)|\bigr), \hspace{0.15cm} \forall x \in \mathbb{R}, 
\end{equation} where $C_0 := 2 + \frac{2}{\phi'(1)} 
\ge 1$. Moreover, $\psi''$ is unbounded. In fact, $\psi''(x_j)=\alpha\hspace{0.07cm}x_j\hspace{0.01cm}\left(1+x_j^2\right)^{(\alpha-2)/2}\sim \alpha \hspace{0.07cm}x_j^{\alpha-1}$ as $j\to\infty$, for $x_j = \left((2\pi j)^{2/\alpha} - 1\right)^{1/2},$ with $j\in\mathbb{N}$.
\end{proposition}

\begin{proof}
We clearly have $\psi \in C^\infty(\mathbb{R})$ with $\psi'(x) = 1 + \sin(\phi(x)) \in [0,2]$. Moreover, \[
\phi'(t) = \alpha\hspace{0.07cm} t\hspace{0.07cm}(1+t^2)^{\alpha/2-1} > 0, \qquad
\phi''(t) = \alpha\hspace{0.07cm}(1+t^2)^{\alpha/2-2}\hspace{0.08cm}\bigl(1 + (\alpha-1)t^2\bigr) > 0,
\] for all $t>0$, because $\alpha > 1$. Hence $\phi'$ is positive and strictly increasing on $[1,\infty[$.  
Set $I(x) := \int_0^x \sin(\phi(t))\,dt$. Let $x\geq 1$. Substituting $u = \phi(t)$ ($du = \phi'(t)\,dt$) and integrating by parts gives
\[
\int_1^x \sin(\phi(t))\,dt
=
\int_{\phi(1)}^{\phi(x)} \frac{\sin u}{\phi'(\phi^{-1}(u))}\,du
=
\left[-\frac{\cos u}{\phi'(\phi^{-1}(u))}\right]_{\phi(1)}^{\phi(x)}
+
\int_{\phi(1)}^{\phi(x)} \cos u \,\frac{d}{du}\!\left(\frac{1}{\phi'(\phi^{-1}(u))}\right)du.
\]
On the one hand, since $u \mapsto 1/\phi'(\phi^{-1}(u))$ is positive on $[\phi(1),\phi(x)]$, we have  \begin{equation*}
    \left|\left[-\frac{\cos u}{\phi'(\phi^{-1}(u))}\right]_{\phi(1)}^{\phi(x)}\right| \leq \frac{1}{\phi'(1)} + \frac{1}{\phi'(x)}.
\end{equation*} 
On the other hand, since $u \mapsto 1/\phi'(\phi^{-1}(u))$ is strictly decreasing on $[\phi(1),\phi(x)]$, we have 
\[
\left|\int_{\phi(1)}^{\phi(x)} \cos u \,\frac{d}{du}\!\left(\frac{1}{\phi'(\phi^{-1}(u))}\right)du\right|\leq \int_{\phi(1)}^{\phi(x)} \left(-\frac{d}{du}\frac{1}{\phi'(\phi^{-1}(u))}\right)du
=
\frac{1}{\phi'(1)} - \frac{1}{\phi'(x)}.
\]
 Therefore, we get $$|I(x)|\leq  \left|\int_0^1 \sin(\phi(t))\,dt\right| + \left|\int_1^x \sin(\phi(t))\,dt\right| \le 1 + \frac{2}{\phi'(1)} =: M,$$ for all $x\geq 1$. Since $g=\sin\circ \phi$ is even, $I(x)= \int_0^x \sin(\phi(t))\,dt$ is odd. Hence, $|I(x)|\leq M,$ for all $|x|\geq 1$. \\

Notice that $M=1+\frac{2}{\alpha \cdot 2^{\frac{\alpha}{2} - 1}}>1$. Clearly, $|I(x)|\leq |x|<1\leq M$ for all $|x|< 1$. Hence $|I(x)|=|\psi(x)-x| \le M$, for all $x\in\mathbb{R}$, so $1+|x| \le 1+M+|\psi(x)| \le C_0(1+|\psi(x)|)$, for all $x\in\mathbb{R}$, where $C_0 := 1+M = 2 + \frac{2}{\phi'(1)}$.\\

An easy computation gives $$\psi''(x) = \alpha \hspace{0.05cm}x\hspace{0.05cm} (1+x^2)^{\frac{\alpha}{2}-1} \hspace{0.05cm}\cos\left((1+x^2)^{\frac{\alpha}{2}}\right),$$ for all $x\in\mathbb{R}$. We now set $x_j := \left((2\pi j)^{2/\alpha} - 1\right)^{1/2},$ for all $j\in\mathbb{N}$. Obviously, we have $\psi''(x_j)= \alpha \hspace{0.07cm}x_j\hspace{0.01cm}\left(1+x_j^2\right)^{(\alpha-2)/2}\sim \alpha \hspace{0.07cm} x_j^{\alpha-1}\to \infty$, as $j\to \infty$, since $\alpha>1$.
\end{proof}

We rely heavily on the following estimate, which can be derived using Cauchy's inequality or Faà di Bruno's formula twice:
\begin{proposition}\label{prop:cauchy}
If $\psi$ is the function given in \eqref{counterexamplePsi}, then there is $B\geq 2$ such that
\begin{equation}\label{eq:psi_deriv}
|\psi^{(\ell)}(x)| \le B^\ell\,\hspace{0.04cm}\ell!\,\hspace{0.04cm}\left(1+|x|\right)^{\beta(\ell-1)},
\end{equation} for all $\ell\in\mathbb{N}$ and $x\in\mathbb{R}$.
\end{proposition}

\begin{proof}
Consider the closed strip $\Omega := \{\zeta \in \mathbb{C} : |\Im \zeta| \le 1/2\}$. For any $\zeta = u + iv \in \Omega$, we have $v^2 \le 1/4$ and $|\zeta|^2 = u^2+v^2 \le u^2 + 1/4$. Therefore,
\begin{equation}\label{AuxComplexI}
\Re(1+\zeta^2) = 1 + u^2 - v^2 \ge \frac{3}{4} + u^2 \ge \frac{1}{2}\left(1 + u^2 + \frac{1}{4}\right) \ge \frac{1}{2}\bigl(1 + |\zeta|^2\bigr) \geq \frac{1}{2}.
\end{equation}
Since $\Re(1+\zeta^2) > 0$ on an open neighborhood of $\Omega$, the principal logarithm $\log(1+\zeta^2)$ (with imaginary part between $-\pi$ and $\pi$) is holomorphic on an open neighborhood of $\Omega$, and thus
\[
\phi(\zeta) := \exp\left(\frac{\alpha}{2}\log(1+\zeta^2)\right), \qquad g(\zeta) := \sin(\phi(\zeta))
\]
are holomorphic on an open neighborhood of $\Omega$ and coincide with $\phi(x)$ and $g(x)$ for $x\in \mathbb{R}$, respectively. \\

Fix $x \in \mathbb{R}$ and define the disk $D_x := \{z \in \mathbb{C} : |z-x| \le \rho_x\}$ centered at $x$ of radius
\[
\rho_x := c\,(1+|x|)^{-\beta},
\] for a suitable $0<c\leq \frac{1}{2}$ to be chosen later. Since $\beta > 0$ and $c \le 1/2$, we have $\rho_x \le 1/2$, and hence $D_x \subset \Omega$. \\

For any $\zeta \in D_x$, note that $\frac{\alpha}{2}-1 < 0$ and $|1+\zeta^2| \ge \Re(1+\zeta^2) \ge \frac{1}{2}(1+|\zeta|^2)$ (see \eqref{AuxComplexI}), which implies the following inequality:
\begin{align*}
|\phi'(\zeta)|
&= \alpha\hspace{0.07cm}|\zeta|\,\exp\!\left(\left(\frac{\alpha}{2}-1\right)\Re\bigl(\log(1+\zeta^2)\bigr)\right) \\
&= \alpha\hspace{0.07cm}|\zeta|\,|1+\zeta^2|^{\alpha/2-1} \\
&\le \alpha\hspace{0.07cm}|\zeta|\left(\frac{1}{2}(1+|\zeta|^2)\right)^{\alpha/2-1} \\
&\le \alpha\hspace{0.07cm} 2^{1-\alpha/2}\,(1+|\zeta|^2)^{1/2}\,(1+|\zeta|^2)^{\alpha/2-1} \\
&\le \alpha\hspace{0.07cm} 2^{1-\alpha/2}\,(1+|\zeta|)^\beta.
\end{align*}
Moreover, $|\zeta| \le |x| + \rho_x \le |x| + \frac{1}{2}$ implies $1+|\zeta| \le \frac{3}{2} + |x| < 2(1+|x|)$. Consequently, for every $\zeta \in D_x$,
\begin{equation}\label{InAuxComplex}
|\phi'(\zeta)|
\le
\alpha\hspace{0.07cm}2^{1-\alpha/2}\,2^{\alpha-1}\,(1+|x|)^\beta
=
\alpha \hspace{0.07cm}2^{\alpha/2}\,(1+|x|)^\beta.
\end{equation}
Clearly, the segment $[x,z]$ lies in $D_x$, for all $z \in D_x$. Since $\phi(x) \in \mathbb{R}$, \eqref{InAuxComplex} implies that
\begin{equation*}
\begin{split}
|\Im\phi(z)| &= |\Im(\phi(z)-\phi(x))| \le \left| \int_{[x,z]} \phi'(\zeta)\,d\zeta \right| \\&\le |z-x|\sup_{\zeta\in D_x}|\phi'(\zeta)| \le \rho_x \cdot \alpha \, 2^{\alpha/2} \, (1+|x|)^\beta \\&= c\,\hspace{0.03cm}\alpha\,\hspace{0.03cm}2^{\alpha/2},
\end{split}
\end{equation*} for all $z \in D_x$. If $c>0$ is chosen to be less than or equal to $\frac{1}{\alpha \hspace{0.07cm} 2^{\alpha/2}}$, then $|\Im\phi(z)|\leq 1$, for all $z\in D_x$, and hence,
\[
|g(z)| = |\sin(\phi(z))| \le e^{|\Im\phi(z)|} \le e, \qquad \forall\hspace{0.05cm}z \in D_x.
\]

Set $c:=\min\{\frac{1}{2}, \frac{1}{\alpha \hspace{0.07cm} 2^{\alpha/2}}\}$ and $B:=\frac{1}{c}=\max\{2, \alpha \hspace{0.07cm} 2^{\alpha/2}\}\geq 2$. By applying Cauchy's inequality to $g$ on $\partial D_x = \{z\in\mathbb{C} : |z-x|=\rho_x\}$ we obtain
\begin{equation}\label{CauchyInMain}
|g^{(n)}(x)|
=
\left|\frac{n!}{2\pi i}\oint_{|z-x|=\rho_x} \frac{g(z)}{(z-x)^{n+1}}\,dz\right|
\le
\frac{e\,n!}{\rho_x^n}
=
e\,B^n\,n!\,(1+|x|)^{\beta n},
\end{equation} for every $n \ge 0$ and $x\in\mathbb{R}$.\\

For $\ell = 1$, we have $|\psi'(x)| = |1 + g(x)| \le 2$, so \eqref{eq:psi_deriv} holds since $B \ge 2$. For $\ell \ge 2$, we take $n=\ell-1 \ge 1$. Clearly, $\psi^{(\ell)}(x)=g^{(n)}(x)$ for all $x\in\mathbb{R}$. By applying \eqref{CauchyInMain}, we obtain:
\[
|\psi^{(\ell)}(x)| \le e\,B^{\ell-1}\,(\ell-1)!\,(1+|x|)^{\beta(\ell-1)} \le B^\ell\,\ell!\,(1+|x|)^{\beta(\ell-1)},
\] for all $x\in\mathbb{R}$. Thus \eqref{eq:psi_deriv} holds for all $\ell \ge 1$ and $x \in \mathbb{R}$.
\end{proof}

From Proposition~\ref{prop:cauchy}, we obtain that \[
 \prod_{\ell=1}^j \frac{\bigl|\psi^{(\ell)}(x)\bigr|^{k_\ell}}{\ell!^{k_\ell}}
\le
\prod_{\ell=1}^j \Bigl(B^\ell\,(1+|x|)^{\beta(\ell-1)}\Bigr)^{k_\ell}
=
B^j\,(1+|x|)^{\beta(j-k)},
\] for all $x\in\mathbb{R}$, $j\in\mathbb{N}$ and $(k_1,...,k_j)\in\mathbb{N}_0^j$ such that $k_1+2k_2+...+j k_j=j$, where $k=k_1+...+k_j$. \\

Using Faà di Bruno's formula and the fact that (see \cite{AFG24}) $$\mathlarger{\sum_{}}\frac{(k_1+...+k_j)!}{k_1!\ldots \hspace{0.05cm} k_j!} = 2^{j-1},$$  where the summation is taken over all $(k_1,...,k_j)\in\mathbb{N}_0^j$ satisfying $k_1+2k_2+\cdots+jk_{j}=j$, 
we get the crucial estimate
\begin{equation}\label{eq:faa_bound}
\bigl|(f\circ\psi)^{(j)}(x)\bigr|
\le
(2B)^j \hspace{0.07cm}\mathlarger{\sum_{k=1}^{j}} \frac{j!}{k!}\,
\bigl|f^{(k)}(\psi(x))\bigr|\,(1+|x|)^{\beta(j-k)},
\end{equation} for all $x\in\mathbb{R}$ and $j\in\mathbb{N}$. \\

We need the following technical lemma:
\begin{lemma}\label{lem:gamma}
Let $d>1$ and $0<\beta\le 1-\frac{1}{d}$. For every $j\in\mathbb{N}$ and $1\leq k\leq j$, set $m := k + \lfloor\beta(j-k)\rfloor + 1$. Then, 
\[
\frac{j!}{k!}\left(\frac{m!}{j!}\right)^d \le (j+1)^d \le 2^{dj}.
\]
\end{lemma}

\begin{proof}

Since $0<\beta<1$ and $1 \le k \le j$, we have $m = k + \lfloor\beta(j-k)\rfloor + 1 \le k + (j-k) + 1 = j + 1$ and also
\[
1 \le k \le m - 1 = k + \lfloor\beta(j-k)\rfloor \le k + \beta(j-k) = (1-\beta)k + \beta j.
\]
It is well known that the function $x \mapsto \log\Gamma(x)$ is strictly increasing on $[2,\infty[$ and convex on $[1,\infty[$ (see for instance \cite{Gamma}). From $\Gamma(m)=(m-1)!$ and $(1-\beta)k + \beta j + 1 = (1-\beta)(k+1) + \beta(j+1)$, we obtain
\begin{equation}\label{logGamma}
\log\bigl((m-1)!\bigr)
\le
\log\Gamma\bigl((1-\beta)(k+1) + \beta(j+1)\bigr)
\le
(1-\beta)\log(k!) + \beta\log(j!).
\end{equation}
By multiplying by $d$ and taking exponentials in \eqref{logGamma}, we obtain $$\bigl((m-1)!\bigr)^d \le (k!)^{d(1-\beta)}\,(j!)^{d\beta},$$ and hence,  
\[
\frac{j!}{k!}\left(\frac{m!}{j!}\right)^d
=
m^d\,\hspace{0.05cm}\frac{j!}{k!}\left(\frac{(m-1)!}{j!}\right)^d
\le
m^d\,\hspace{0.04cm}\frac{(k!)^{d(1-\beta)}\,(j!)^{d\beta}}{k!\,(j!)^{d-1}}
=
m^d\hspace{0.07cm}\left(\frac{k!}{j!}\right)^{d(1-\beta)-1}.
\]
Since $1 \le k \le j$, we have $\frac{k!}{j!} \le 1$. Moreover, the hypothesis $\beta \le 1 - \frac{1}{d}$ is precisely the condition $d(1-\beta) - 1 \ge 0$, which ensures that $\bigl(\frac{k!}{j!}\bigr)^{d(1-\beta)-1} \le 1$. Therefore,
\[
\frac{j!}{k!}\left(\frac{m!}{j!}\right)^d
\le
m^d
\le
(j+1)^d
\le
2^{dj}, 
\] where we have used that $j+1\leq 2^j$ in the last inequality.

\end{proof}

\begin{proof}[Proof of Theorem~\ref{thm:main}] Part~1 was already proved in Proposition~\ref{prop:growth}. \\

To prove Part~2, let $f \in \Sigma_d(\mathbb{R})$ and fix $h > 0$. We will find $T > 0$ (depending only on $h, d, \alpha$) such that $\|C_\psi f\|_{h,d} \le 2^d\hspace{0.09cm} \|f\|_{T,d}$.\\

First, for $j=0$ and $Q \in \mathbb{N}_0$, \eqref{eq:comparable} and the definition of $\|f\|_{T,d}$ (see \eqref{eq:seminorm}) give
\begin{equation}\label{eq:j0}
(1+|x|)^Q\,|f(\psi(x))|
\le
C_0^Q\,(1+|\psi(x)|)^Q\,|f(\psi(x))|
\le
(C_0T)^Q\,Q!^{\,d}\,\|f\|_{T,d},
\end{equation} for all $x\in\mathbb{R}$. \\

Now let $j \ge 1$ and $Q \in \mathbb{N}_0$. Multiplying \eqref{eq:faa_bound} by $(1+|x|)^Q$ gives
\begin{equation}\label{mainInequalityOsc}
(1+|x|)^Q\,\bigl|(f\circ\psi)^{(j)}(x)\bigr|
\le
(2B)^j\hspace{0.07cm} \mathlarger{\sum_{k=1}^{j}} \frac{j!}{k!}\hspace{0.03cm}\,(1+|x|)^{Q+\beta(j-k)}\hspace{0.03cm}\,\bigl|f^{(k)}(\psi(x))\bigr|,
\end{equation} for all $x\in\mathbb{R}$. \\

For each $1\leq k \leq j$, we define
\[
p_k := \lfloor\beta(j-k)\rfloor + 1, \qquad m\equiv m(k) := k + p_k \le j + 1, \qquad q\equiv q(k) := Q + p_k.
\]

Since $Q + \beta(j-k) < q\leq  m+Q$, $1+|x| \ge 1$, and $C_0 \ge 1$, \eqref{eq:comparable} implies that $$(1+|x|)^{Q+\beta(j-k)} \le C_0^q \hspace{0.09cm}(1+|\psi(x)|)^q \le C_0^{m+Q}\hspace{0.09cm}(1+|\psi(x)|)^q,$$ for all $x\in\mathbb{R}$. \\

Since $k + q = m + Q$, the following inequality holds:
\begin{equation}\label{key1Osc}
(1+|x|)^{Q+\beta(j-k)}\,\hspace{0.02cm}\bigl|f^{(k)}(\psi(x))\bigr|
\le
\|f\|_{T,d}\hspace{0.02cm}\,(C_0T)^{m+Q}\,\hspace{0.02cm}(m+Q)!^{\,d},
\end{equation} for all $x\in\mathbb{R}$. \\

Now, using $(m+Q)! \le 2^{m+Q}\,\hspace{0.02cm}m!\,\hspace{0.02cm}Q! \le 2^{j+1+Q}\,\hspace{0.02cm}m!\,\hspace{0.02cm}Q!$, Lemma~\ref{lem:gamma}, and $j!\,\hspace{0.02cm}Q! \le (j+Q)!$, we obtain the following inequality: 
\begin{equation}\label{Key2Osc}
\frac{j!}{k!}\,(m+Q)!^{\,d}
\le
2^{d(j+1+Q)}\left(\frac{j!}{k!}(m!)^d\right)(Q!)^d
\le
2^d \cdot 4^{dj} \cdot 2^{dQ}\,(j+Q)!^{\,d}.
\end{equation}

If $C_0T \le 1$, then the inequality $m \ge k + \beta(j-k) = (1-\beta)k + \beta j > \beta j$ (as $\beta < 1$) implies $m+Q > \beta(j+Q)$, and therefore $(C_0T)^{m+Q} \le \bigl((C_0T)^\beta\bigr)^{j+Q}$.\\

By inserting \eqref{key1Osc} and \eqref{Key2Osc} into \eqref{mainInequalityOsc}, we finally get:
\begin{equation}\label{eq:jge1}
\begin{split}
(1+|x|)^Q\,\bigl|(f\circ\psi)^{(j)}(x)\bigr|
&\le
(2B)^j \sum_{k=1}^{j} 2^d \cdot 4^{dj} \cdot 2^{dQ}\,\Bigl((C_0T)^\beta\Bigr)^{j+Q} (j+Q)!^{\,d}\,\|f\|_{T,d}
\\ &\le
2^d\,\Bigl(4B\cdot 4^d\,(C_0T)^\beta\Bigr)^j\,\Bigl(2^d\,(C_0T)^\beta\Bigr)^Q\,(j+Q)!^{\,d}\,\|f\|_{T,d},
\end{split}
\end{equation} for all $x\in\mathbb{R}$, where we have used that $\sum_{k=1}^j 1=j\leq 2^j$. \\

We choose
\[
T := \frac{1}{C_0} \min\!\left\{1,\; h,\; \left(\frac{h}{4B\cdot 4^d}\right)^{1/\beta}\right\} > 0.
\]
Then $C_0T \le 1$, $C_0T \le h$ (which yields $(C_0 T)^Q \le h^Q$ in \eqref{eq:j0} for $j=0$), and $4B\cdot 4^d\,\hspace{0.02cm}(C_0T)^\beta \le h$, which implies $2^d\hspace{0.07cm}(C_0T)^\beta \le h$ since $B \ge 2$. Notice that the constant $T>0$ does not depend on $f$. \\

Therefore, combining \eqref{eq:j0} and \eqref{eq:jge1} and taking the supremum over $x \in \mathbb{R}$ and $j, Q \in \mathbb{N}_0$, we conclude that
\[
\|C_\psi f\|_{h,d} \le 2^d\,\|f\|_{T,d} < \infty.
\]

Since $h>0$ was arbitrary, $C_\psi : \Sigma_d(\mathbb{R}) \to \Sigma_d(\mathbb{R})$ is continuous.
\end{proof}

For the symbol $\psi$ defined in \eqref{counterexamplePsi} with $1 < \alpha \le 2 - \frac{1}{d}$ and $\beta = \alpha - 1 \in \left]0, 1 - \frac{1}{d}\right]$, Proposition~\ref{prop:cauchy} shows that $|\psi^{(m)}(x)| \le B^m\hspace{0.05cm} m!\hspace{0.05cm} (1 + |x|)^{\beta(m-1)}$ for every $m \ge 2$ and $x \in \R$. We also verified that the second derivative $\psi''$ attains its maximum polynomial growth of degree $\beta$. For $\beta=1-\frac{1}{d}$, in Theorem \ref{mainth} we proved that this growth coincides with the largest allowed for symbols in $\Sigma_d(\mathbb{R})$. We now show that the same sharpness holds for every higher-order derivative.
\begin{proposition}\label{remarkpsider}
Let $d > 1$, $1 < \alpha \le 2 - \frac{1}{d}$, and $\beta := \alpha - 1 \in \left]0, 1 - \frac{1}{d}\right]$. For the symbol $\psi$ defined in \eqref{counterexamplePsi} and every integer $m \ge 2$, there exists a sequence $(x_j^{(m)})_{j \in \N}$ tending to $+\infty$ such that
\[
\bigl|\psi^{(m)}(x_j^{(m)})\bigr| \sim \alpha^{m-1}\,(x_j^{(m)})^{\beta(m-1)} \qquad \text{as } j \to \infty.
\]
In particular, when $\alpha = 2 - \frac{1}{d}$, $\psi^{(m)}$ attains polynomial growth of degree $(m-1)\bigl(1 - \frac{1}{d}\bigr)$.
\end{proposition}
\begin{proof}
Fix $m \ge 2$ and set $n := m - 1 \ge 1$, so that $\psi^{(m)} = g^{(n)} = (\sin \circ \phi)^{(n)}$. First, we bound the derivatives of $\phi(x) = (1+x^2)^{\alpha/2}$. For each $x \in \R$, consider the closed disk
\[
\widetilde{D}_x := \left\{ \zeta \in \C : |\zeta - x| \le \frac{1}{4}(1 + |x|) \right\}.
\] 
Let $\zeta = u + iv \in \{\eta \in \C : |\eta - x| \le \frac{1}{3}(1 + |x|)\}$. If $|x| \le 1$, then $|v| \le \frac{2}{3}$ and hence $\Re(1 + \zeta^2) = 1 + u^2 - v^2 \ge \frac{5}{9} > 0$. If $|x| > 1$, then $|u| \ge |x| - \frac{1}{3}(1+|x|) = \frac{2|x|-1}{3} > 0$ and $|v| \le \frac{1}{3}(1+|x|)$, which implies
\[
\Re(1+\zeta^2) \ge 1 + \frac{(2\vert{}x\vert{}-1)^2 - (1+\vert{}x\vert{})^2}{9} = \frac{3\vert{}x\vert{}^2 - 6\vert{}x\vert{} + 9}{9} = \frac{1}{3}(\vert{}x\vert{}-1)^2 + \frac{2}{3} > 0.
\]
 Therefore, $\phi(\zeta) := \exp\bigl(\frac{\alpha}{2}\log(1+\zeta^2)\bigr)$ is holomorphic on an open neighborhood of $\widetilde{D}_x$ and coincides with $\phi(x)$ for $x\in\mathbb{R}$. Moreover, for every $\zeta \in \widetilde{D}_x$ we have $|\zeta| \le |x| + \frac{1}{4}(1+|x|) \le \frac{5}{4}(1+|x|)$, and thus
\[
|\phi(\zeta)| \le (1 + |\zeta|^2)^{\alpha/2} \le \left(1 + \frac{25}{16}(1+|x|)^2\right)^{\alpha/2} \le 3^{\alpha/2}\hspace{0.07cm}(1+|x|)^\alpha.
\]
By applying Cauchy's inequality to $\phi$ on $\widetilde{D}_x$ we obtain:
\begin{equation}\label{eq:phi_cauchy_remark}
|\phi^{(\ell)}(x)| \le \frac{\ell!\, 3^{\alpha/2}\hspace{0.07cm}(1+|x|)^\alpha}{\bigl(\frac{1}{4}(1+|x|)\bigr)^\ell} = 3^{\alpha/2}\, \hspace{0.03cm}4^\ell\,\hspace{0.03cm} \ell!\,\hspace{0.03cm} (1+|x|)^{\alpha - \ell}, 
\end{equation} for every $\ell \in \N$ and $x \in \R$. \\

Now we apply Fa\`a di Bruno's formula to $\psi^{(m)}(x) = (\sin \circ \phi)^{(n)}(x)$ and separate the unique term corresponding to $k = n$ (for which $k_1 = n$ and $k_2 = \dots = k_n = 0$):
\[
\psi^{(m)}(x) = \sin^{(n)}(\phi(x))\, \bigl(\phi'(x)\bigr)^n + R_n(x),
\]
where $R_1(x) \equiv 0$ when $m = 2$, and for $m \ge 3$,
\[
R_n(x) := \sum_{k=1}^{n-1} \sin^{(k)}(\phi(x)) \sum_{\substack{k_1+\dots+k_n = k \\ k_1+2k_2+\dots+nk_n = n}} \frac{n!}{k_1!\cdots k_n!} \prod_{\ell=1}^{n} \left(\frac{\phi^{(\ell)}(x)}{\ell!}\right)^{k_\ell}.
\]
For every $1 \le k \le n - 1$ and every $(k_1,\dots,k_n)\in\mathbb{N}_0^{n}$ with $\sum_{\ell=1}^n k_\ell = k$ and $\sum_{\ell=1}^n \ell k_\ell = n$, \eqref{eq:phi_cauchy_remark} implies
\[
\prod_{\ell=1}^{n} \left(\frac{|\phi^{(\ell)}(x)|}{\ell!}\right)^{k_\ell} \le \prod_{\ell=1}^{n} \Bigl(3^{\alpha/2}\, 4^\ell\, (1+|x|)^{\alpha - \ell}\Bigr)^{k_\ell} = 3^{\alpha k / 2}\, 4^n\, (1+|x|)^{\alpha k - n}.
\]
Since $k \le n - 1$ and $\alpha > 1$, we have $\alpha\hspace{0.05cm} k - n \le \alpha\hspace{0.05cm}(n-1) - n = (\alpha - 1)\hspace{0.05cm}n - \alpha = \beta\hspace{0.05cm}(m-1) - \alpha$. As $|\sin^{(k)}(\phi(x))| \le 1$, there exists a constant $K_{m,\alpha} \ge 0$ (depending only on $m$ and $\alpha$) such that
\begin{equation}\label{eq:rem_bound_remark}
|R_n(x)| \le K_{m,\alpha}\,(1+|x|)^{\beta(m-1) - \alpha} \qquad \forall x \in \R.
\end{equation}

Finally, we define the sequence $(x_j^{(m)})_{j \in \N}$ by
\[
x_j^{(m)} := \left((2\pi j + \theta_m)^{2/\alpha} - 1\right)^{1/2}, \qquad \theta_m := \begin{cases} 0 & \text{if } m \text{ is even}, \\ \pi/2 & \text{if } m \text{ is odd}. \end{cases}
\]
Because $2\pi j + \theta_m \ge 2\pi > 1$, each $x_j^{(m)} > 0$ is well-defined, $x_j^{(m)} \to +\infty$ as $j \to \infty$, and $\phi(x_j^{(m)}) = 2\pi j + \theta_m$ (for $m=2$, $x_j^{(2)}$ coincides with $x_j$ in Theorem~\ref{thm:main}(1)). Since $\sin^{(n)}(u) = \sin\bigl(u + \frac{n\pi}{2}\bigr)$ with $n = m - 1$, we have $|\sin^{(n)}(2\pi j)| = |\cos(2\pi j)| = 1$ when $m$ is even ($n$ odd) and $|\sin^{(n)}(2\pi j + \pi/2)| = |\sin(2\pi j + \pi/2)| = 1$ when $m$ is odd ($n$ even). Consequently, $|\sin^{(n)}(\phi(x_j^{(m)}))| = 1$ for all $j \in \N$. Moreover, $\phi'(x) = \alpha \hspace{0.07cm}x (1+x^2)^{(\alpha-2)/2} \sim \alpha \hspace{0.07cm}x^{\alpha - 1} = \alpha\hspace{0.07cm} x^\beta$ as $x \to +\infty$, so
\begin{equation} \label{AsympHO}
    \bigl|\sin^{(n)}(\phi(x_j^{(m)}))\,\bigl(\phi'(x_j^{(m)})\bigr)^n\bigr| = \bigl(\phi'(x_j^{(m)})\bigr)^{m-1} \sim \alpha^{m-1}\,(x_j^{(m)})^{\beta(m-1)} \qquad (j \to \infty).
\end{equation}
Combining \eqref{AsympHO} with \eqref{eq:rem_bound_remark} and using $\beta(m-1) - \alpha < \beta(m-1)$ (since $\alpha > 1$), we conclude that
\[
\bigl|\psi^{(m)}(x_j^{(m)})\bigr| \sim \alpha^{m-1}\,(x_j^{(m)})^{\beta(m-1)} \qquad \text{as } j \to \infty.
\]
When $\alpha = 2 - \frac{1}{d}$, $\psi^{(m)}$ attains the polynomial growth of degree $(m-1)\bigl(1 - \frac{1}{d}\bigr)$.
\end{proof}

\begin{remark}
In light of these results, we confirm that the sufficient condition in \cite[Corollary 4.2]{AFG24} for the composition operator $C_\psi: \Sigma_d(\mathbb{R})\to \Sigma_d(\mathbb{R})$ to be continuous (which in particular requires $|x| \le C(1 + |\psi(x)|)$ and all derivatives of $\psi$ to be bounded) is far from necessary, leaving substantial room to relax its hypotheses.   
\end{remark}
We expect that an extension of the techniques in Section~\ref{mainsec} yields the following conjecture: \begin{conjecture}
Let $d > 1$ and suppose $C_\psi : \Sigma_d(\R) \to \Sigma_d(\R)$ is well-defined and continuous. Then for every $m\in\mathbb{N}$ and $\varepsilon > 0$, there exists a constant $C_{m,\varepsilon} > 0$ such that
\begin{equation*}
|\psi^{(m)}(y)| \le C_{m,\varepsilon}\hspace{0.05cm} (1 + |y|)^{(m-1)\hspace{0.03cm}(1 - \frac{1}{d}) + \varepsilon} \qquad \forall y \in \R.
\end{equation*}    
\end{conjecture}

\end{document}